\documentclass[12t, russian]{article}
\usepackage[cp1251]{inputenc}
\usepackage[T2A]{fontenc}
\usepackage[english,russian]{babel}
\usepackage{amsmath, amssymb}
\usepackage{multirow}
\usepackage{xcolor}
\usepackage{graphicx}

\newtheorem{theorem}{Theorem}
\newtheorem{lemma}{Lemma}
\newtheorem{proposition}{Proposition}

\voffset=- 30mm
\begin{document}
\sloppy
\large


\begin{center}
\textbf{Differentiation of resultants and common roots of pairs of polynomials}
\end{center}
\bigskip
\bigskip

\begin{center}
\textbf{\large M. M. Chernyavsky (1), A. V. Lebedev (2),  Yu. V. Trubnikov (1)}
\end{center}

\bigskip

\begin{center}
{\large ((1) Vitebsk State University,  Belarus, \	(2) Belarusian State University,  Belarus)}
\end{center}

\bigskip
\bigskip
\bigskip

\begin{abstract}

{
The well-known mathematical instrument for detection  common roots for pairs of polynomials
and multiple roots of polynomials are resultants and discriminants.
For a pair of polynomials   $f$ and $g$ their resultant $R(f,g)$ is a function of their coefficients.  Zeros of resultant  $R(f,g)$ correspond
   to the families of coefficients of  $f$ and $g$ such that  $f$ and $g$ have a common root.
Herewith the calculation of this common root is a separate problem.

The principal results on calculation of  a unique common root of two polynomials and
also about calculating a unique root of multiplicity 2 of a polynomial
in terms of the first order partial derivatives of resultants and discriminants
are given in  the monograph by I.M.~Gelfand, \ M.M.~Kapranov, \
A.V.~Zelevinsky~\cite[Ch.~3, Ch.~12]{GKZ} 
  A significant development of the ideas
of this book in the direction of searching for formulas for multiple roots of
polynomials is presented in the paper by  I. A. Antipova,  E. N. Mikhalkin,  A. K. Tsikh~\cite{AMTs}. The key result of this article is \cite[Theorem~1]{AMTs}
where the expression for
a unique root of multiplicity $s \geq 3$ in terms of the first order partial derivatives
of resultant of the polynomial and it's derivative of order  $s-1$.  

In the present article the explicit formulas for  higher derivatives of resultants
of pairs of polynomials possessing common roots are obtained. On this basis a series of
 results that differ in ideas from  \cite[Theorem~1]{AMTs} linking  higher derivatives of resultants and common multiple  roots  are proven.
In addition the results obtained are applied for a new transparent proof of
a refinement of \cite[Theorem~1]{AMTs}.	
}

\end{abstract}
\bigskip
\bigskip


\bigskip
\bigskip


\bigskip
\bigskip
\bigskip


\medskip

\medskip

\section*{Introduction}

The well-known mathematical instrument for detection  common roots for pairs of polynomials
and multiple roots of polynomials are resultants and discriminants.
For a pair of polynomials
  $f$ and $g$ their resultant $R(f,g)$ is a function of their coefficients
  (for precise definition see section~\ref{s-1}). Zeros of resultant  $R(f,g)$ correspond
   to the families of coefficients of  $f$ and $g$ such that  $f$ and $g$ have a common root.
Herewith the calculation of this common root is a separate problem.
Up to the end of the XX century finding the explicit formulas expressing the values of
common (multiple) roots of polynomials in terms of their coefficients has been a hard task since
in general it needs voluminous analytical calculation not amenable to manual counting.
A relatively new approach in this direction is expressing the value of a common root as well as
a multiple root (so also certain combinations of multiple roots if a multiple root is not unique)
in terms of partial derivatives of some resultants and discriminants.
Among the authoritative sources in the theory of algebraic equations that also contains
information about calculating a unique common root of two polynomials and
also about calculating a unique root of multiplicity 2 of a polynomial
in terms of the first order partial derivatives of resultants and discriminants
one should mention  the monograph by I.M.~Gelfand, \ M.M.~Kapranov, \
A.V.~Zelevinsky~\cite[Ch.~3, Ch.~12]{GKZ} (see Proposition~\ref{p-GKZ}
and Theorem~\ref{t-5} of the present article).
Herewith the authors of~\cite{GKZ} emphasize that in their opinion
'The general proof given
there (reduction via the Cayley trick, to the biduality theorem for projective dual
varieties), is probably the most transparent'.
 (see  \cite[c. 400]{GKZ}).  A significant development of the ideas
of this book in the direction of searching for formulas for multiple roots of
polynomials is presented in the paper by  I. A. Antipova,  E. N. Mikhalkin,  A. K. Tsikh~\cite{AMTs}. The key result of this article is \cite[Theorem~1]{AMTs}
(Theorem~\ref{t-6} of the present paper) where the expression for
a unique root of multiplicity $s \geq 3$ in terms of the first order partial derivatives
of resultant of the polynomial and it's derivative of order  $s-1$.  The proof of Theorem~\ref{t-6}
in this article uses a series of specific methods of algebra and complex analysis
 including well-known formulas for resultants, residue theory and the form of polynomial
 $f$ in Hilbert strata. This technique, of course, corresponds to the fundamental nature of the topic under consideration (cf. the quoted remark of the authors of~\cite{GKZ}), but it makes it difficult for applied specialists to perceive the material of the work.

Note also that  current algorithms for machine calculation
of multiple roots of polynomials considered, for example, in~\cite{Ch} and~\cite{YaCh},
in essence do not give any novelty from the point of view of
representing  multiple roots in an explicit form in terms of coefficients.

In the present article the explicit formulas for  higher derivatives of resultants
of pairs of polynomials possessing common roots are obtained. On this basis a series of
 results that differ in ideas from  Theorem~\ref{t-6} linking  higher derivatives of resultants and common multiple  roots  are proven.
In addition the results obtained are applied for a new transparent proof of
a refinement of Theorem~\ref{t-6}.	

By concrete examples (in Section~\ref{s-3}) we illustrate
qualitative differences between explicit formulas for multiple roots of polynomials
obtained on the basis of the results of \cite{GKZ,AMTs} and the results of this article.

\section{Resultants and common roots of pairs of polynomials}\label{s-1}

\subsection{Resultant and discriminant}\label{ss-1}

Here we recall the definitions and properties of the resultant and discriminant necessary for further presentation. For details see, for example,  \cite[$\S$ 54]{Kur}.	
	
Let  $f\left(z \right):=\sum\limits_{i=0}^n {a_i {z^{n-i}}}\, \left(a_0\ne 0 \right)$
	and $g \left(z \right):=\sum\limits_{j=0}^m {b_j {z^{m-j}}}\,\, \left(b_0 \ne 0 \right)$
be two polynomials having  roots   $\left\{\alpha _i \right\}_{i=1}^n$ and $\left\{\beta _j \right\}_{j=1}^m$ respectively. \textit{Resultant} of the polynomials  $f$ and $g$ is the product , \begin{equation}\label{e-0}
 R(f,g ):=a_0^mb_0^n\cdot \prod\limits_{i=1}^n{\prod\limits_{j=1}^m{\left({\alpha _i}-{\beta _j} \right)}}.
\end{equation}

Polynomials  $f$ and  $g$ in the definition of the resultant
are used in the definition of  resultant in a non-symmetrical way. Clearly
\begin{equation}\label{e-0**}
 R(g ,f)=  (-1)^{mn}R(f,g) .
\end{equation}

If the polynomials $f$ and $g$ possess at least one common root then
$R(g ,f)= 0$. This property of resultant can also be used to check for multiple roots of a polynomial, to do this, it is enough to calculate the resultant of the polynomial and its first derivative.

Since  $f\left(z \right)=a_0 \prod\limits_{i=1}^n{\left(z-{\alpha _i} \right)},\quad g\left(z \right)=b_0 \prod\limits_{j=1}^m{\left(z-{\beta _j} \right)}$ their resultant
 $R(f,g)$ can be represented in the form
\begin{equation}\label{e-1}
R(f,g) =a_0^m  g \left(\alpha _1 \right)g \left(\alpha _2 \right)\cdot \ldots \cdot g \left(\alpha _n \right)={{\left(-1 \right)}^{mn}}b_0^n f\left(\beta _1 \right)f\left(\beta _2 \right)\cdot \ldots \cdot f\left(\beta _m \right).
\end{equation}

It is also known that the resultant can be represented as the following determinant, which is sometimes called  \textit{Sylvester formula} or  \textit{resultant in Sylvester form}
	\begin{equation}\label{e-S}
		R(f,g) =\left| \begin{matrix}
		\medskip
		 a_0 & a_1 & \ldots  & a_n & {}  \\
{}&a_0 & a_1 & \ldots  & a_n & {} \\
\cdots  & \cdots  & \cdots  & \cdots  & \cdots  & \cdots   \\		
				\medskip
		{} & {} & a_0 & a_1 & \ldots  & a_n  \\
		\medskip
		b_0 & b_1 & \ldots  & b_m & {} & {}  \\
		\medskip
		{} & b_0 & b_1 & \ldots  & b_m & {}  \\
		
		\cdots  & \cdots  & \cdots  & \cdots  & \cdots  & \cdots   \\
		\medskip
		{} & {} & b_0 & b_1 & \ldots  & b_m  \\
		
	\end{matrix}
	\right| \ .
	\end{equation}
In this determinant $m$ lines containing coefficients of polynomial $f\left(z \right)$
and  $n$ lines containing coefficients of polynomial $g\left(z \right)$.
In the cells of the determinant, where empty fields are left, zeros are implied.

\smallskip

Formula  \eqref{e-S} naturally leads us to consider the resultant as
a function of the coefficients of the polynomials
 $f$ and $g$, i.e.  $R(f,g)(a_0, ..., a_n, b_0, ..., b_m)$.
Namely  in this sense  we will understand the resultant in this work.

\smallskip	
	
 For  a polynomial   $
f\left(z \right)=a_0 z^n+a_1 {z^{n-1}}+\ldots +{a_n}\quad \left(a_0 \ne 0 \right)$,
 having roots  $\left\{\alpha _i \right\}_{i=1}^n$ its \emph{discriminant} $D(f)$ is defined as
$$
D(f) := a_0^{2n-2}\prod_{i>j}(\alpha_i -\alpha_j)^2\,.
$$
The value of the resultant of a polynomial and its first derivative
is proportional to the value of its discriminant  $D(f)$. Namely, the following relation
is true
\begin{equation}\label{e-*}	
R\left(f, f' \right)={{\left(-1 \right)}^{n\left(n-1 \right)/2}}a_0 \,D(f).
\end{equation}

\bigskip

\subsection{Resultants and common multiple roots of polynomials}\label{ss-2}

In the monograph by I.M.~Gelfand, \ M.M.~Kapranov, \
A.V.~Zelevinsky~\cite{GKZ} the following statement is proven  (it follows from Corollary 3.7  \cite[с.~109]{GKZ}).

\begin{proposition}\label{p-GKZ}
If polynomials
		\begin{equation}\label{e-15}
f\left( x \right)=a_0{x^n}+a_1{x^{n-1}}+a_2{x^{n-2}}+\ldots +a_{n-1}x+a_n
		\end{equation}		
		and
\begin{equation}\label{e-16}
		g\left( x \right)=b_0x^m+b_1x^{m-1}+b_2x^{m-2}+\ldots +b_{m-1}x+b_m
\end{equation}
have exactly only one
simple common root $w$, then  $w$  can be found explicitly from relations$:$
\begin{equation}\label{e-17}
		\left( {{w}^n}: {{w }^{n-1}}  :{{w }^{n-2}} : \ldots  \ :  w   :1 \right)=\left( \frac{\partial R\left( f,g \right)}{\partial a_0}  : \frac{\partial R\left( f,g \right)}{\partial a_1} : \ldots :  \frac{\partial R\left( f,g \right)}{\partial {a_n}} \right),
\end{equation}
\begin{equation}\label{e-18}
		\left( {{w }^m} : {{w }^{m-1}}: {{w }^{m-2}} :\ldots  : w : 1 \right)=\left( \frac{\partial R\left( f,g \right)}{\partial b_0} : \frac{\partial R\left( f,g \right)}{\partial {b_1}}: \ldots : \frac{\partial R\left( f,g \right)}{\partial {b_m}} \right),
\end{equation}		
where  $R\left( f,g \right)$ is the resultant of the polynomials $f\left( x \right)$ and $g\left( x \right)$.
\end{proposition}	

In  \cite{GKZ} formals  \eqref{e-17}, \eqref{e-18} are derived by
interpretation of discriminant in the duality theory of projective
varieties, and the authors of~\cite{GKZ} emphasize that in their opinion this proof is  the most transparent (see discussion  \cite[p.~400]{GKZ}).

To satisfy relations \eqref{e-17} and \eqref{e-18}
it is fundamentally important that there is only one common root,
and that it is not a multiple for any of the polynomials.
Here one can give an example  $f\left( x \right)={{\left( x-1 \right)}^3},\quad g\left( x \right)={{\left( x-1 \right)}^2}$, where all the first order partial derivatives  $\frac{\partial R\left( f,g \right)}{\partial {a_i}},\, i=0,1,2,3; \ \frac{\partial R\left( f,g \right)}{\partial {b_j}},\, j=0,1, 2;$ are zero.

\smallskip

This section is devoted to describing the effects that occur in such situations.

\medskip

Let us consider polynomials
 $f\left(z \right) =\sum\limits_{i=0}^n {a_i {z^{n-i}}}$ $\left(a_0\neq 0\right)$
 and  $g \left(z \right):=\sum\limits_{j=0}^m {b_j {z^{m-j}}}\,\, \left(b_0 \ne 0 \right)$.

 In accordance with  \eqref{e-1} we write the resultant  $R\left(f , g \right)$
 in the form

\begin{equation}\label{e-3-}
	 	R\left(f,g \right)=:R=a_0^mg_1 g_2\cdot \ldots \cdot {g_{n-1}} g_n,
\end{equation}
where
\begin{equation}\label{e-4-}
	g_i \equiv g\left(z_i \right)= b_0 z_i^{m}+ b_1 z_i^{m-1}+\ldots +{b_{m-1}}z_i+{b_{m}}=\sum\limits_{j=0}^{m}{b_j}z_i^{m-j}\quad \left(i=1, 2,\ldots , n \right)
\end{equation}	
is the value of polynomial   $g$ at  \textit{i}-th root of polynomial  $f$ \
($z_i, \ i=1, ... , n$ are the roots of  $f$).

In addition, again in accordance with  \eqref{e-1}
\begin{equation}\label{e-3--}
R(f,g) ={{\left(-1 \right)}^{mn}}b_0^n f_1\cdot\ldots\cdot f_m,
\end{equation}
where
\begin{equation}\label{e-4--}
	f_i \equiv f\left(y_i \right)= a_0 y_i^{n}+ a_1 y_i^{n-1}+\ldots +{a_{n-1}}y_i+{a_{n}}=\sum\limits_{j=0}^{n}{a_j}y_i^{n-j}\quad \left(i=1, 2,\ldots , m \right)
\end{equation}	
is the value of polynomial  $f$ at \textit{i}-th root of polynomial $g$ \
($y_i, \ i=1, ... , m$ -- корни $g$).

\medskip

Our first observation is

\begin{lemma}\label{l*-1}
Let  $z_1=w$ be a root of a polynomial  $f\left(z \right)=\sum\limits_{i=0}^n {a_i {z^{n-i}}}$ $\left(a_0\neq 0 \right)$ \
 $($$z_i, \ i=1, ..., n$ are the roots of   $f$$)$
 and let $w=y_1$ be also a root  of a polynomial   $g \left(z \right):=\sum\limits_{j=0}^m {b_j {z^{m-j}}}\,\, \left(b_0 \ne 0 \right)$ \
 $($$y_i, \ i=1, ... , m$ are the roots of  $g$$)$.
 Then for the resultant
 $R:= R\left(f,g \right)$ the following equalities are true
  \begin{equation}\label{e-l*-1}
 \frac{\partial R}{\partial {b_j}} = a_0^m w^{m-j}  \prod\limits_{i=2}^n g_i, \ \ \ j=0, ... , m;
 \end{equation}
 where  $g_i$ are give by  \eqref{e-4-}$;$

  \begin{equation}\label{e-l*-1*}
 \frac{\partial R}{\partial {a_j}} = {{\left(-1 \right)}^{mn}}b_0^n w^{n-j}  \prod\limits_{i=2}^m f_i, \ \ \ j=0, ... , n;
 \end{equation}
 where $f_i$ are given by  \eqref{e-4--}.
 \end{lemma}

 P\,r\,o\,o\,f. \ In accordance with  \eqref{e-3-}

\begin{equation}\label{e-*1-l}
 R = a_0^m g_1 \cdot  \prod\limits_{i=2}^n g_i\,.
\end{equation}
Since  $z_1=w$ then according to \eqref{e-4-} one has
\begin{equation}\label{e-4*-l}
	g_1 =  \sum\limits_{j=0}^{m}{b_j}w^{m-j}.
\end{equation}
And therefore
\begin{equation}\label{e-4**-l}
\frac{\partial g_1}{\partial b_{j}} = w^{m-j}.
\end{equation}
Thus
$$
\frac{\partial R}{\partial {b_j}} = a_0^m\frac{\partial }{\partial {b_j}}\left( g_1 \cdot  \prod\limits_{i=2}^n g_i\right) = a_0^m \frac{\partial g_1}{\partial {b_j}} \prod\limits_{i=2}^n g_i + a_0^m g_1 \frac{\partial }{\partial {b_j}}\left(\prod\limits_{i=2}^n g_i\right)=
   a_0^m w^{m-j}  \prod\limits_{i=2}^n g_i,
$$
where in the final equality we equality  we take into account
 \eqref{e-4**-l} and the fact that  $g_1 =g(w)=0$. The equality \eqref{e-l*-1} is proven.

The equality  \eqref{e-l*-1*} can be proven by the same reasoning
using the equalities
\eqref{e-3--}  and~\eqref{e-4--}. The proof is finished.

\medskip

Lemma~\ref{l*-1} implies the next refinement of Proposition~\ref{p-GKZ}.

\begin{proposition} \label{p-2} \
Let  $f\left(z \right)=\sum\limits_{i=0}^n {a_i {z^{n-i}}}$ $\left(a_0\neq 0, \ a_n \neq 0   \right)$,
and   $g \left(z \right):=\sum\limits_{j=0}^m {b_j {z^{m-j}}}\,\, \left(b_0 \neq 0, \  \ b_m \neq 0\right)$.
There is a unique simple common root $w$ is  of the polynomials  $f$ and $g$
iff the following two conditions are satisfied$:$

$1)$   $R(f,g)=0;$

$2)$
$ \frac{\partial  R}{\partial b_{m}}\neq 0 $, \  $\frac{\partial  R}{\partial a_{n}}\neq 0$.

\smallskip

Herewith  $w$ satisfies relations  \eqref{e-17} and \eqref	{e-18}.
\end{proposition}

P\,r\,o\,o\,f. \  Necessity. Let $w=z_1=y_1$ be a common root of  $f$ and $g$.
Thus  $R(f,g)=0$. Since  $a_n \neq 0, \ b_m \neq 0$ all the roots of the polynomials
 $f$ and $g$ are non-zero. Let  $z_1, \ldots , z_n$ be the roots of  $f$, and
 $y_1, \ldots , y_m$ be the roots of $g$. Since
$z_2, \ldots , z_n$ are not the roots of  $g$ it follows that
\begin{equation}\label{e-*1-l-}
  \prod\limits_{i=2}^n g_i \neq 0\,,
\end{equation}
and as
$y_2, \ldots , y_m$ are not the roots of $f$ we have that
\begin{equation}\label{e-l*-1*-}
\prod\limits_{i=2}^m f_i \neq 0\,.
 \end{equation}
From inequalities  \eqref{e-*1-l-} and  \eqref{e-l*-1*-} taking into account
\eqref{e-l*-1} and \eqref{e-l*-1*} we have
$$
\frac{\partial  R}{\partial b_{j}}\neq 0, \quad j=1,\ldots, m; \qquad\frac{\partial  R}{\partial a_{j}}\neq 0, \quad j=1,\ldots, n;
$$
i.e. condition 2) is satisfied.

Sufficiency.
 Since  $R(f,g)=0$ the polynomials $f$ and $g$ have a common root, say,
$z_1=y_1=w$. Since  $\frac{\partial  R}{\partial b_{m}}\neq 0$ the formula  \eqref{e-l*-1}
implies
$$
\prod\limits_{i=2}^n g_i \neq 0 .
$$
Thus  $z_2, \ldots , z_n$ are not the roots of  polynomial  $g$ (in particular,
 $z_i\neq w=z_1, \ i=2, \ldots , n$) and $w$ ia a simple root of $f$.

Similarly since  $\frac{\partial  R}{\partial a_{n}}\neq 0$ the formula
 \eqref{e-l*-1*} implies
$$
\prod\limits_{i=2}^m f_i \neq 0 .
$$
Therefore   $y_2, \ldots , y_m$ are not the roots of  polynomial $f$  and $w$
ia a simple root of $g$.

Finally, relations  \eqref{e-17} and \eqref	{e-18} follow directly from
 \eqref{e-l*-1*} and \eqref{e-l*-1}. The proof is complete.

\medskip

The next statement describes higher derivatives of resultants of polynomials possessing
common multiple roots.
	
\begin{theorem}\label{t**-1}
Let  $z_1=z_2= ... z_s=w$ be a root of multiplicity at least~$s$ $(2  \leq s<n)$
of a polynomial
$f\left(z \right)=\sum\limits_{i=0}^n {a_i {z^{n-i}}}$
$\left(a_0\neq 0 \right)$ \
 $($here  $z_i, \ i=1, ..., n$ are the roots of $f$$)$.
 Let  $w$ be also a root of a polynomial
  $g \left(z \right):=\sum\limits_{j=0}^m {b_j {z^{m-j}}}\,\, \left(b_0 \ne 0 \right)$.
Then for the resultant
 $R:= R\left(f,g \right)$ the following equalities are true

\medskip

 $1)$ for  $1 \leq r< s$
 \begin{equation}\label{e-*0-}
  \frac{\partial ^r R}{\partial b_{j_r}...\partial b_{j_1}}=0, \ \ \ \ \ \ \ \ \left(j_k=0, 1, \ldots , m \right);
  \end{equation}

 $2)$
\begin{equation}\label{e-*00-}
 \frac{\partial ^s R}{\partial b_{j_s}...\partial b_{j_1}}=a_0^ms! {w^{sm-\left(j_s+...+j_1 \right)}}\cdot \prod\limits_{i=s+1}^n g_i\ \  \ \ \ \ \ \ \left(j_k=0, 1, \ldots , m \right),
\end{equation}
where $g_i$ are given by \eqref{e-4-}.

\smallskip

If  $s=n$, i.e. $f(z) =a_0 (z-w)^n$, then for $r<n$ equality \eqref{e-*0-}
takes place and
\begin{equation}\label{e-*000-}
 \frac{\partial ^n R}{\partial b_{j_n}...\partial b_{j_1}}=a_0^m n! {w^{nm-\left(j_n+...+j_1 \right)}} \ \  \ \ \ \ \ \ \left(j_k=0, 1, \ldots , m \right).
\end{equation}
\end{theorem}

P\,r\,o\,o\,f. \ In accordance with  \eqref{e-3-} we have

\begin{equation}\label{e-*1-}
 R = a_0^m g_1 \cdot ...\cdot g_s \cdot \prod\limits_{i=s+1}^n g_i\,.
\end{equation}
Since  $z_1=z_2= ... z_s=w$ it follows by  \eqref{e-4-} that
\begin{equation}\label{e-4*-}
	g_i =  \sum\limits_{j=0}^{m}{b_j}w^{m-j}=: h \quad \ \ \ \ \ \left(i=1, 2,\ldots , s \right).
\end{equation}
Thus, in particular,
\begin{equation}\label{e-4**-}
\frac{\partial h}{\partial b_{j}} = w^{m-j}.
\end{equation}
Taking into account  \eqref{e-4*-} and denoting  $\prod\limits_{i=s+1}^n g_i =: P$
one can write  \eqref{e-*1-} in the form
\begin{equation}\label{e-*2-}
 R = a_0^m h^s \cdot P\,.
\end{equation}

Let us verify the equality  \eqref{e-*0-}.

 \begin{equation*}
  \frac{\partial ^r R}{\partial b_{j_r}...\partial b_{j_1}}=  a_0^m \frac{\partial ^r }{\partial b_{j_r}...\partial b_{j_1}} (h^s\cdot P) =
    a_0^m \frac{\partial ^{r-1} }{\partial b_{j_r}...\partial b_{j_2}} \left(\frac{\partial}{\partial b_{j_1}}(h^s\cdot P)\right) =
  \end{equation*}
  $$
=  a_0^m \frac{\partial ^{r-1} }{\partial b_{j_r}...\partial b_{j_2}} \left(sh^{s-1}\cdot
\frac{\partial h}{\partial b_{j_1}}\cdot P + h^s\frac{\partial P}{\partial b_{j_1}}\right) =
  $$
\begin{equation}\label{e-3*-}
=  a_0^m \frac{\partial ^{r-1} }{\partial b_{j_r}...\partial b_{j_2}} \left(w^{m-j_1}sh^{s-1}
\cdot P + h^s\frac{\partial P}{\partial b_{j_1}}\right) ,
\end{equation}
where in the final equality we used  \eqref{e-4**-}.

From \eqref{e-3*-} by continuing implementing differentiation
one obtains
\begin{equation*}
 \frac{\partial ^r R}{\partial b_{j_r}...\partial b_{j_1}}=
  a_0^m \frac{\partial ^{r-1} }{\partial b_{j_r}...\partial b_{j_2}} \left(w^{m-j_1}sh^{s-1}
\cdot P + h^s\frac{\partial P}{\partial b_{j_1}}\right)  =
\end{equation*}
$$
= a_0^m \frac{\partial ^{r-2} }{\partial b_{j_r}...\partial b_{j_3}} \left(\frac{\partial}{\partial b_{j_2}}\left(w^{m-j_1}sh^{s-1}
\cdot P + h^s\frac{\partial P}{\partial b_{j_1}}\right) \right) =
$$
$$
= a_0^m \frac{\partial ^{r-2} }{\partial b_{j_r}...\partial b_{j_3}} \left(
\left(w^{m-j_1}s(s-1)h^{s-2}\frac{\partial h}{\partial b_{j_2}}
\cdot P +  w^{m-j_1}sh^{s-1}
\cdot \frac{\partial P}{\partial b_{j_2}}\right)+ \frac{\partial}{\partial b_{j_2}}\left(h^s\frac{\partial P}{\partial b_{j_1}}\right) \right) =
$$
\begin{equation}\label{e-*5-}
= a_0^m \frac{\partial ^{r-2} }{\partial b_{j_r}...\partial b_{j_3}} \left(
\left(w^{2m-(j_1+j_2)}s(s-1)h^{s-2}
\cdot P +  w^{m-j_1}sh^{s-1}
\cdot \frac{\partial P}{\partial b_{j_2}}\right)+ \frac{\partial}{\partial b_{j_2}}\left(h^s\frac{\partial P}{\partial b_{j_1}}\right) \right),
\end{equation}
where in the final equality we again used \eqref{e-4**-}.

The equality \eqref{e-*5-} is rewritten in the form

\begin{equation*}
 \frac{\partial ^r R}{\partial b_{j_r}...\partial b_{j_1}}= a_0^m \left( w^{2m-(j_1+j_2)}s(s-1)\frac{\partial ^{r-2} }{\partial b_{j_r}...\partial b_{j_3}} \left(h^{s-2}\cdot P\right) + \right.
\end{equation*}
\begin{equation}\label{e-*6-}
 \left. +  w^{m-j_1}s\frac{\partial ^{r-2} }{\partial b_{j_r}...\partial b_{j_3}}\left(h^{s-1}
\cdot \frac{\partial P}{\partial b_{j_2}}\right)+ \frac{\partial ^{r-1} }{\partial b_{j_r}...\partial b_{j_2}}
\left(h^s\frac{\partial P}{\partial b_{j_1}}\right)\right).
\end{equation}
Note that (since $r<s$) in any summand in \eqref{e-*6-}
the order of differentiation is less than the degree of the multiplier $h$
included in this term. And herewith the procedure of differentiation of every
summand is the same as the procedure of differentiation of the initial function:
\begin{equation*}
    \frac{\partial ^r }{\partial b_{j_r}...\partial b_{j_1}} (h^s\cdot P).
\end{equation*}
Continuing the procedure of differentiation  \eqref{e-*6-}, i.e. implementing
$r$ steps of differentiation and taking into account the foregoing observation we arrive
at the expression
\begin{equation}\label{e-*7-}
 \frac{\partial ^r R}{\partial b_{j_r}...\partial b_{j_1}}= a_0^m w^{rm-(j_1+...+j_r)}s(s-1)\cdot ...\cdot (s-r+1) \left(h^{s-r}\cdot P\right) + \dots ,
\end{equation}
where by  $\dots$ we denoted the summands each of which
  (taking into account the above observation)
contains as a multiplier $h^k, \ k>0$.

Since  $w$ is a common root of the polynomials $f$ and $g$
it follows in accordance with \eqref{e-4-} and
\eqref{e-4*-}, that \mbox{$h=g(w)=0$}. Thus  \eqref{e-*7-} implies
\begin{equation*}
 \frac{\partial ^r R}{\partial b_{j_r}...\partial b_{j_1}}= 0.
\end{equation*}
The equality \eqref{e-*0-} is proven.

\smallskip

Let us verify the equality  \eqref{e-*00-}.

For  $r=s$ the equality  \eqref{e-*7-} turns into

\begin{equation}\label{e-*8-}
 \frac{\partial ^s R}{\partial b_{j_s}...\partial b_{j_1}}= a_0^m w^{sm-(j_1+...+j_s)}s! \cdot P + \dots ,
\end{equation}
where, as above, by  $\dots$ we denoted the summands each of which
  (taking into account the above observation)
contains as a multiplier $h^k, \ k>0$.

Again taking into account  $h=g(w)=0$
we obtain \eqref{e-*00-}.

\smallskip

Note finally that for  $s=n$
the equality  \eqref{e-*2-} turns into the equality
\begin{equation*}
 R = a_0^m h^n ,
\end{equation*}
that directly implies  \eqref{e-*000-}.

The proof is complete.

\medskip

If in Theorem~\ref{t**-1}  we replace the polynomials $f$ and $g$ in places
then taking into account  \eqref{e-3--}
and \eqref{e-4--} one obtains the next statement.

\begin{proposition}\label{p-4}
Let  $y_1=y_2= ... y_s=w$ be a root of multiplicity at least~$p$ $(2  \leq p<m)$
for a polynomial
$g \left(z \right):=\sum\limits_{j=0}^m {b_j {z^{m-j}}}\,\, \left(b_0 \ne 0 \right)$
 $($here  $y_j, \ j=1, ..., m$ are the roots of  $g$$)$.
Let $w$ be also a root of a polynomial
  $f\left(z \right)=\sum\limits_{i=0}^n {a_i {z^{n-i}}}$ $\left(a_0\neq 0 \right)$.
Then for the resultant
 $R:= R\left(f,g \right)$ the following equalities are true

 \smallskip

 $1)$ for  $1 \leq r< p$
 \begin{equation}\label{e-**0-}
  \frac{\partial ^r R}{\partial a_{i_r}...\partial a_{i_1}}=0, \ \ \ \ \ \ \ \ \left(i_k=0, 1, \ldots , n \right);
  \end{equation}

$2)$
\begin{equation}\label{e-*00---*}
 \frac{\partial ^p R}{\partial a_{i_p}...\partial a_{i_1}}={{\left(-1 \right)}^{mn}}b_0^np! {w^{pn-\left(i_p+...+i_1 \right)}}\cdot \prod\limits_{j=p+1}^m f_j\ \  \ \ \ \ \ \ \left(i_k=0, 1, \ldots , n \right),
\end{equation}
where  $f_j$ are given by  \eqref{e-4--}.

\smallskip

If  $p=m$, i.e.     $g(z) =b_0 (z-w)^m$, then for  $r<m$ equality  \eqref{e-**0-}
takes place and

\begin{equation}\label{e-**000-}
 \frac{\partial ^m R}{\partial a_{i_m}...\partial a_{i_1}}= b_0^n m! {w^{nm-\left(i_m+...+i_1 \right)}} \ \  \ \ \ \ \ \ \left(i_k=0, 1, \ldots , n \right).
\end{equation}
\end{proposition}

As a corollary of the forgoing observations one can at once obtain a generalization of
Proposition~\ref{p-GKZ}
on the case of common multiple roots of two polynomials.

\begin{theorem}\label{t-r1}
Let $f\left(z \right)=\sum\limits_{i=0}^n {a_i {z^{n-i}}}$ $\left(a_0\neq 0, \ a_n \neq 0   \right)$, and
   $g \left(z \right):=\sum\limits_{j=0}^m {b_j {z^{m-j}}}\,\, \left(b_0 \neq 0, \  \ b_m \neq 0\right)$.
If the polynomials  $f$ and $g$ have a single common root  $w$ of multiplicity  $s$ for  $f$
and multiplicity  $p$
for  $g$
then  $w$ satisfies the relations$:$
\begin{equation}\label{e-*00--}
 \frac{\partial ^s R}{\partial b_{j_s}...\partial b_{j_1}} : \frac{\partial ^s R}{\partial b_{k_s}...\partial b_{k_1}} = {w^{\left(k_s+...+k_1 \right)-\left(j_s+...+j_1 \right)}} \qquad\left(j_r, k_r=0, 1, \ldots , m\right).
\end{equation}
and
\begin{equation}\label{e-*00---}
 \frac{\partial ^p R}{\partial a_{i_p}...\partial a_{i_1}} : \frac{\partial ^p R}{\partial a_{l_p}...\partial a_{l_1}} = {w^{\left(l_p+...+l_1 \right)-\left(i_p+...+i_1 \right)}} \qquad\left(i_r, l_r=0, 1, \ldots , n\right).
\end{equation}
\end{theorem}

P\,r\,o\,o\,f. \   Relations  \eqref{e-*00--} follow directly from
 \eqref{e-*00-} taking into account that in the situation under consideration  $\prod\limits_{i=s+1}^n g_i \neq 0$; and relations  \eqref{e-*00---} follow from
\eqref{e-*00---*} taking into account that in the situation under consideration
$\prod\limits_{j=p+1}^m f_j \neq 0$.
The proof is complete.

\medskip

\subsection{Differential calculus of resultants and
common roots of pairs of polynomials}\label{ss-1a}

The results of the previous section show that
 in the presence of common roots for pairs of polynomials,
 the derivatives of their resultants take on a rather special form.
In this section we will write these results in the corresponding terms.

Recall that for a function  $F: {\mathbb{C}}^p \to \mathbb{C} \ (p=1,2, \dots )$
its derivative of order  $s$ $F^{(s)}$ \ $(s=1,2, \dots)$ is an $s-$linear function
given by the family of partial derivatives
$$
F^{(s)} \leftrightarrow \left[\frac{\partial ^s F}{\partial z_{i_s}...\partial z_{i_1}}\right], \ \ i_r= 1,\dots , p\,.
$$
If  $F(a,b), \ a=(a_0, \dots, a_n), \  b= (b_0, ..., b_m)$ is a function of
 $(n+1)+(m+1)$ variables, then by $F^{(s)}_a$ we denote it's partial derivative
 of order $s$ on variables~$a$, and, respectively, by  $F^{(s)}_b$
  we denote it's partial derivative
 of order  $s$ on variables~$b$. For the first partial derivatives
 on $a$ and $b$ it is natural also to use the notation
$\frac{\partial  F}{\partial a}$ and $\frac{\partial  F}{\partial b}$.

Let as write the results obtained in section~\ref{ss-2} in this terms.
\smallskip

Lemma~\ref{l*-1} takes the following form.

\begin{lemma}\label{*l*-1}
Let  $z_1=w$ be a root of a polynomial
$f\left(z \right)=\sum\limits_{i=0}^n {a_i {z^{n-i}}}$ $\left(a_0\neq 0 \right)$ \
 $($$z_i, \ i=1, ..., n$ are the roots of    $f$$)$, and let $w=y_1$
 be also a root  of a polynomial
  $g \left(z \right):=\sum\limits_{j=0}^m {b_j {z^{m-j}}}\,\, \left(b_0 \ne 0 \right)$ \ $($$y_i, \ i=1, ... , m$ are the roots of $g$$)$.  Then  for the resultant
 $R:= R\left(f,g \right)\,(a,b)$, \ $a= (a_0, \dots, a_n)$,  \ $b=(b_0, \dots, b_m)$
 the following equalities are true
  \begin{equation}\label{e-*l*-1}
 \frac{\partial R}{\partial {b}} = a_0^m  \prod\limits_{i=2}^n g_i \cdot [w^{m-j}], \ \ \ j=0, ... , m;
 \end{equation}
 where $g_i$ are given by \eqref{e-4-}$;$

  \begin{equation}\label{e-*l*-1*}
 \frac{\partial R}{\partial {a}} = {{\left(-1 \right)}^{mn}}b_0^n   \prod\limits_{i=2}^m f_i \cdot [w^{n-j}], \ \ \ j=0, ... , n;
 \end{equation}
 where $f_i$ are given by \eqref{e-4--}.
 \end{lemma}

Theorem~\ref{t**-1} is written in the following form.

\begin{theorem}\label{*t**-1}
Let $z_1=z_2= ... z_s=w$ be a root of multiplicity at least~$s$ $(2  \leq s<n)$
for a polynomial
$f\left(z \right)=\sum\limits_{i=0}^n {a_i {z^{n-i}}}$ $\left(a_0\neq 0 \right)$ \
 $($here  $z_i, \ i=1, ..., n$ are the roots of  $f$$)$.
 Let  $w$ be also a root of a polynomial   $g \left(z \right):=\sum\limits_{j=0}^m {b_j {z^{m-j}}}\,\, \left(b_0 \ne 0 \right)$.  Then for the resultant
 $R:= R\left(f,g \right) (a,b)$ \ $a= (a_0, \dots, a_n)$,  \ $b=(b_0, \dots, b_m)$
the following equalities are true

\medskip

 $1)$ for  $1 \leq r< s$
 \begin{equation}\label{*e-**0-}
   R^{(r)}_b=0,
  \end{equation}

 $2)$
\begin{equation}\label{*e-**00-}
  R^{(s)}_b =a_0^ms!  \prod\limits_{i=s+1}^n g_i  \cdot\left[{w^{sm-\left(j_s+...+j_1 \right)}}\right],
\end{equation}
where  $g_i$ are given by  \eqref{e-4-}.

\smallskip

If  $s=n$, i.e.    $f(z) =a_0 (z-w)^n$, then for  $r<n$ equality  \eqref{*e-**0-}
takes place, and
\begin{equation}\label{*e-**000-}
 R^{(n)}_b =a_0^m n! \cdot \left[{w^{nm-\left(j_n+...+j_1 \right)}}\right]\,.
\end{equation}
\end{theorem}

Proposition~\ref{p-4} is written in the following form respectively.

\begin{proposition}\label{*p-4}
Let  $y_1=y_2= ... y_s=w$ be a root of multiplicity at least~$p$ $(2  \leq p<m)$
for a polynomial
$g \left(z \right):=\sum\limits_{j=0}^m {b_j {z^{m-j}}}\,\, \left(b_0 \ne 0 \right)$
 $($here  $y_j, \ j=1, ..., m$ are the roots of  $g$$)$.
Let $w$ be also a root of a polynomial
  $f\left(z \right)=\sum\limits_{i=0}^n {a_i {z^{n-i}}}$ $\left(a_0\neq 0 \right)$.
Then for the resultant
 $R:= R\left(f,g \right)$ the following equalities are true

 \smallskip

 $1)$ for  $1 \leq r< p$
 \begin{equation}\label{**e-**0-}
   R^{(r)}_a=0 ,
  \end{equation}

$2)$
\begin{equation}\label{**e-*00---*}
 R^{(p)}_a ={{\left(-1 \right)}^{mn}}b_0^np!  \prod\limits_{j=p+1}^n f_j \cdot \left[ {w^{pn-\left(i_p+...+i_1 \right)}}\right],
\end{equation}
where $f_j$ are given by  \eqref{e-4--}.

\smallskip

If $p=m$, i.e.     $g(z) =b_0 (z-w)^m$, then for  $r<m$ the equality \eqref{**e-**0-}
takes place and

\begin{equation}\label{**e-**000-}
 R^{(m)}_a = b_0^n m! \cdot\left[{w^{nm-\left(i_m+...+i_1 \right)}}\right] .
\end{equation}
\end{proposition}

\section{Rational expressions for multiple roots of polynomials}\label{s-2}

\subsection{Algorithms of derivation of explicit rational expressions for
multiple roots of polynomials (known results)}	\label{ss-3}

	Here we recall the well-known algorithms of derivation of explicit rational expressions for
multiple roots of polynomials.

Since a multiple root of a polynomial is a root of its derivative the relations
\eqref{e-17} and \eqref{e-18} can be used for deriving
the formulas for this root calculation obtaining in this way series
of formulas  ``in the spirit of Viete's formulas''.
In~\cite{GKZ} we have the next statement.

\begin{theorem}\label{t-5}  {\rm{\cite[с.~404]{GKZ}}}.
If for a polynomial
 $f\left( x \right)={a_0}{x^{n}}+{a_1}{x^{n-1}}+{a_2}{x^{n-2}}+\ldots +{a_n}$
 its discriminant $D\left( f \right)$ is zero but at least
one partial derivative   $\partial D\left( f \right)/\partial {a_i}\ \ \left( i=0,\ 1,\ \ldots ,\ n \right)$ is non-zero, then  $f$ has a unique  root  $w$ of multiplicity 2,
and it can be found from the proportions$:$
\begin{equation}\label{e-19}
		\left( {{w }^n} : {{w }^{n-1}}: \ldots  : w  : 1 \right)=\left( \frac{\partial D\left( f \right)}{\partial a_0} : \frac{\partial D\left( f \right)}{\partial a_1} : \ldots  : \frac{\partial D\left( f \right)}{\partial a_{n-1}} : \frac{\partial D\left( f \right)}{\partial a_n} \right).
\end{equation}
\end{theorem}
	
Recalling  relation \eqref{e-*} between resultant and discriminant and the fact that
the resultant is 	symmetric (up to a constant factor) with respect to the
polynomials of which it is composed, we note that expressions \eqref{e-17} and \eqref{e-18}
are equal (i.e. one can equally use partial derivatives on coefficients of $f$
and on coefficients on $f'$). The authors in~\cite{GKZ} do not
develop further the idea of calculating of a root of multiplicity 2 by means
of partial derivatives of resultant on coefficients of polynomial $g\left( x \right)$
as, for example, in~\eqref{e-18}.

In the case when the original polynomial $f\left( x \right)$ has at least one multiple root, all partial derivatives of the first order of the resultant
$R\left( f,g \right)$ by the coefficients of the polynomial
$g\left( x \right)={f}'\left( x \right)$
 are equal to zero (Theorem~\ref{t*-1} of the present article).
 In this case, to calculate the value of a multiple root, it is necessary to involve
 higher derivatives (see theorem~\ref{t*-1} of this article).
  In the monograph~\cite{GKZ} no prerequisites for this are given.

The article by I.A. Antipova,  E.N. Mikhalkin and  A.K. Tsikh~\cite{AMTs}
 is devoted to the development of ideas from the book~\cite{GKZ} on the explicit expression of multiple roots of a polynomial. The key result here is
		
\begin{theorem}\label{t-6}	
	{\rm \cite[Теорема~1]{AMTs}}.
If  $s\ge 3$, then a non-zero solution to an
algebraic equation with complex coefficients
\begin{equation}\label{e-20}
f=f\left( z \right)=z^n+a_1z^{n-1}+a_2{z^{n-2}}+\ldots +a_{n-1}z+{a_n}=0
\end{equation}
$z_1$ of multiplicity~$s$ is expressed in terms of partial derivatives
 $\dfrac{\partial R}{\partial {a_i}}$, where
  $R=R\left( f,{f^{\left( s-1 \right)}} \right)$
is the resultant of the polynomial  and its derivative of order~$s-1$
according to the formulas
\begin{equation}\label{e-21}
		{z_1}=\frac{\partial R}{\partial a_{n-(s-2)}}:\frac{\partial R}{\partial a_{n-(s-3)}}=\frac{\partial R}{\partial a_{n-(s-3)}}:\frac{\partial R}{\partial a_{n-(s-4)}} \ldots =\frac{\partial R}{\partial a_{n-2}}:\frac{\partial R}{\partial a_{n-1}}=\frac{\partial R}{\partial a_{n-1}}:\frac{\partial R}{\partial a_n}.
\end{equation}

If  $s=2$, then the non-zero root of multiplicity 2 to the \eqref{e-20} $z_1$
is expressed in terms of partial derivatives  $\dfrac{\partial D}{\partial a_i}$
of the discriminant  $D$ according to the formulas
\begin{equation}\label{e-22}
	{{z}_{1}}=\frac{\partial D}{\partial a_1}:\frac{\partial D}{\partial a_2}= \ldots  =\frac{\partial D}{\partial a_{n-2}}:\frac{\partial D}{\partial a_{n-1}}=\frac{\partial D}{\partial a_{n-1}}:\frac{\partial D}{\partial a_n}\, .
\end{equation}
\end{theorem}
	
The second part of this theorem coincides with Theorem~\ref{t-5}  and expression  \eqref{e-19}.

The first comparison of the expressions \eqref{e-21} and \eqref{e-22} shows that they differ
in the number of relations that give the value of a multiple  root $z_1$. Indeed,
the number of such relations in \eqref{e-21} depends only on the degree of multiplicity
of the root $z_1$ and is equal to $s-2$. There is no such restriction in the expressions \eqref{e-22}.

Note that due to the fact that in the formulation of theorem~\ref{t-6} the polynomial $f\left(z\right)$ is written with the highest coefficient $a_0=1$,
the equality $z_1=\dfrac{\partial D}{\partial a_0}:\dfrac{\partial D}{\partial a_1}$ is lost
in the expressions \eqref{e-22}.
The validity of this equality is proved in the theorem ~\ref{t-5}.

In the article~\cite{AMTs} the proof of  Theorem~\ref{t-6}
 uses a number of specific methods of algebra and complex analysis,
 including well-known formulas for resultants, residue theory and
 the form of the $f$  polynomial on Hilbert strata.
 Such a technique, of course,
 corresponds to the fundamental nature of the topic under consideration,
 but makes it difficult for applied specialists to perceive the material of the work.

In the present paper we prove a series of results that differ in ideas from
Theorem~\ref{t-6} and link higher derivatives of resultants and multiple roots.
In addition these results are applied to obtain a new transparent proof
of a refinement of  Theorem~\ref{t-6}.	

\medskip

\subsection{Multiple roots of a polynomial and resultants}	\label{ss-4}

Here we gather a series of corollaries of the results of  Section~\ref{ss-2}
for the resultants of a polynomial and its derivatives.

Let us consider a polynomial  $f\left(z \right) =\sum\limits_{i=0}^n {a_i {z^{n-i}}}$ $\left(a_0\neq 0\right)$.
In accordance with   \eqref{e-1} we write the resultant of the polynomial
 $f\left(z \right)$ and its  $k$-th  derivative  $f^{(k)}\left(z \right)= \sum\limits_{j=0}^{n-k}{b_j {z^{n-k-j}}}$ in the form:
	
\begin{equation}\label{e-3}
	 	R\left(f,f^{(k)} \right)=:R=a_0^m g_1 g_2\cdot \ldots \cdot {g_{n-1}} g_n =
{{\left(-1 \right)}^{mn}}b_0^n f_1f_2\cdot \ldots \cdot f_{n-k},
\end{equation}
where
\begin{equation}\label{e-4}
	g_i \equiv f^{(k)}\left(z_i \right)= \sum\limits_{j=0}^{n-k}{b_j}z_i^{n-k-j}\quad \left(i=1, 2,\ldots , n \right)
\end{equation}	
is the value of the $k$-th derivative  $f^{(k)}$  of  $f$ at its  \textit{i}-th root \ ($z_i, i=1, ... , n$ are the roots of $f$);
\begin{equation}\label{e-44}
	f_i \equiv f\left(y_i \right)= \sum\limits_{j=0}^{n}{a_j}y_i^{n-j}\quad \left(i=1, 2,\ldots , n-k \right)
\end{equation}
is the value of the polynomial  $f$ at  \textit{i}-th root of its  $k$-th derivative
 $f^{(k)}$\ ($y_i, i=1, ... , n-k$ are the roots of $f^{(k)}$).

\smallskip

An immediate corollary of Lemma~\ref{l*-1} is

\begin{lemma}\label{l**-1}
Let  $z_1=z_2= ... z_s=w$ be a root of multiplicity at least $s$ $(2  \leq s<n)$ for a polynomial
 $f\left(z \right)=\sum\limits_{i=0}^n {a_i {z^{n-i}}}$ $\left(a_0\neq 0 \right)$. Then for the resultant
 $R:= R\left(f,\ f^{(k)} \right)$, where  $f^{(k)}\left(z \right)= \sum\limits_{j=0}^{n-k}{{b_j}\,{z^{n-k-j}}}$, $1 \leq k < s $
  the following equalities are true

 \begin{equation}\label{e-l**-1*}
 \frac{\partial R}{\partial {a_j}} = {{\left(-1 \right)}^{(n-k)n}}b_0^n w^{n-j}  \prod\limits_{i=2}^{n-k} f_i, \ \ \ j=0, ... , n;
 \end{equation}
where $f_i$ are given by  \eqref{e-44}.
 \end{lemma}

An immediate corollary of Theorem~\ref{t**-1} is
	
\begin{theorem}\label{t*-1}
Let  $z_1=z_2= ... z_s=w$ be a root of multiplicity at least $s$ $(2  \leq s<n)$
for a polynomial  $f\left(z \right)=\sum\limits_{i=0}^n {a_i {z^{n-i}}}$ $\left(a_0\neq 0 \right)$. Then for the resultant
 $R:= R\left(f,\ f^{(k)} \right)$, where  $f^{(k)}\left(z \right)= \sum\limits_{j=0}^{n-k}{{b_j}\,{z^{n-k-j}}}$, $1 \leq k < s $
 the following equalities are true

 $1)$ for  $1 \leq r< s$
 \begin{equation}\label{e-*0}
  \frac{\partial ^r R}{\partial b_{j_r}...\partial b_{j_1}}=0, \ \ \ \ \ \ \ \ \left(j_k=0, 1, \ldots , n-k \right);
  \end{equation}

 $2)$
\begin{equation}\label{e-*00}
 \frac{\partial ^s R}{\partial b_{j_s}...\partial b_{j_1}}=a_0^{n-k} s! {w^{s\left(n-k \right)-\left(j_s+...+j_1 \right)}}\cdot \prod\limits_{i=s+1}^n g_i\ \  \ \ \ \ \ \ \left(j_k=0, 1, \ldots , n-k \right),
\end{equation}
where  $g_i$ are given by  \eqref{e-4}.

\smallskip

If  $s=n$, i.e.     $f(z) =a_0 (z-w)^n$ then for  $r<n$  equality  \eqref{e-*0}
is true and
\begin{equation}\label{e-*000}
 \frac{\partial ^n R}{\partial b_{j_n}...\partial b_{j_1}}=a_0^{n-k} n! {w^{n\left(n-k \right)-\left(j_n+...+j_1 \right)}} \ \  \ \ \ \ \ \ \left(j_k=0, 1, \ldots , n-k \right).
\end{equation}
\end{theorem}

\subsection{Applications of the main theorems}\label{ss-5}

Here we gather a series of corollaries of the results of the previous section.

\medskip

Note first that Theorem~\ref{t-5} is a corollary of Lemma~\ref{l**-1}.
It follows from formula \eqref{e-l**-1*} for  $s=2,\ k=1$.
Namely, the following more explicit form of Theorem~\ref{t-5} takes place.

\begin{proposition}\label{p-3}
A polynomial
$f\left(z \right)=\sum\limits_{i=0}^n {a_i{z^{n-i}}}\  \left(a_0\neq 0, \ a_n\neq 0 \right)$
can be represented in the form
\begin{equation}\label{e-8-p}
	f\left(z \right)=a_0{{\left(z-w \right)}^2}\left(z-{z_3} \right)\left(z-{z_4} \right)\ldots \left(z-{z_{n-1}} \right)\left(z-{z_n} \right),
\end{equation}
where  $z_i\neq z_j, \ i\neq j$ и $w\neq z_i;$
	
\noindent iff the  following two relations are satisfied$:$

$1)$ $R\left(f,f' \right)=0$ $($or $D\left(f \right)=0$ in view of \eqref{e-*}$);$

$2)$ for $R=R\left(f,f' \right)$ at least one of the derivatives  $\frac{\partial R}{\partial {a_j}}$ is non-zero.

Herewith the vector of derivatives  $\left[\frac{\partial R}{\partial {a_0}},\frac{\partial R}{\partial {a_1}},
\ldots , \frac{\partial R}{\partial {a_n}}\right]$ has the form
\begin{equation}\label{e-l**-1**p-}
\left[\frac{\partial R}{\partial {a_0}},\frac{\partial R}{\partial {a_1}},
\ldots , \frac{\partial R}{\partial {a_n}}\right] = {{\left(-1 \right)}^{(n-1)n}}(na_0)^n   \prod\limits_{i=2}^{n-1} f_i\,\cdot \,[w^n,w^{n-1}, ... , w, 1]
\end{equation}
and thus relation \eqref{e-19} is true.
\end{proposition}

P\,r\,o\,o\,f.\,
Existence of a multiple root for a polynomial $f$ is equivalent to the condition
$R\left(f,f' \right)=0$ (or $D\left(f \right)=0$ in view of \eqref{e-*}).
Under the satisfaction of this condition let $z_1=z_2 =w$ be a multiple root of $f$.
Since
 $a_n\neq 0$ all the roots of $f$ are non-zero (in particular, $w\neq 0$).
 In accordance with formula  \eqref{e-l**-1*} (for  $s=2, k=1$) we have
 \begin{equation}\label{e-l**-1**}
 \frac{\partial R}{\partial {a_j}} = {{\left(-1 \right)}^{(n-1)n}}b_0^n w^{n-j}  \prod\limits_{i=2}^{n-1} f_i, \ \ \ j=0, ... , n.
 \end{equation}
Since  $b_0 =na_0\neq 0$ it follows that
\begin{equation*}
\frac{\partial R}{\partial {a_j}}  \neq 0 \ \ \Leftrightarrow \ \
\prod\limits_{i=2}^{n-1} f_i \neq 0.
\end{equation*}
Therefore
\begin{equation*}
\frac{\partial R}{\partial {a_j}}  \neq 0 \ \ \Leftrightarrow \ \ f_i= f(y_i) \neq 0, \qquad i=2,\ldots , n-1.
\end{equation*}
Thus a common root of the polynomial  $f$ and its derivative $f'$
is only the root
$y_1=w=z_1=z_2$. Therefore $f$ does not possess any other multiple roots, i.e. all the roots
 $z_3, \ldots , z_n$  are simple.
Herewith  \eqref{e-l**-1**p-} is simply a vector record of  \eqref{e-l**-1**}.

Finally the relation  \eqref{e-19} (taking into account  \eqref{e-*})
follows directly from \eqref{e-l**-1**p-}.
The proof is complete.

\medskip

The next statement follows from Theorem	~\ref{t*-1} and is a certain alternative to
Theorem~\ref{t-5} (Proposition~\ref{p-3}).

\begin{theorem}\label{t-2}
A polynomial  $f\left(z \right)=\sum\limits_{i=0}^n {a_i{z^{n-i}}}\  \left(a_0\neq 0, \ a_n\neq 0 \right)$ can be represented in the form
\begin{equation}\label{e-8}
	f\left(z \right)=a_0{{\left(z-w \right)}^2}\left(z-{z_3} \right)\left(z-{z_4} \right)\ldots \left(z-{z_{n-1}} \right)\left(z-{z_n} \right),
\end{equation}
where $z_i\neq z_j, \ i\neq j$ and  $w\neq z_i;$
	
\noindent	iff  the  following two relations are satisfied$:$
	
\smallskip

	$1)$ $R\left(f,f' \right)=0$ $($либо $D\left(f \right)=0$ в силу \eqref{e-*}$);$
	
	$2)$ $\dfrac{\partial ^2 R\left(f,f' \right)}{\partial b_{n-1}^2}\ne 0,$ where ${b_j}\
\left(j=0, 1,\ldots , n-1 \right)$ are the coefficients of the polynomial $g\left(z \right)=f'\left(z \right)$.

Herewith for the multiple root  $w$ the following relations are satisfied

\begin{equation}\label{e-*0*0}
 \frac{\partial ^2 R}{\partial b_{j_2}\partial b_{j_1}} :
 \frac{\partial ^2 R}{\partial b_{k_2}\partial b_{k_1}} =
 w^{(k_2+k_1) -(j_2 +j_1)}, \qquad
  \left(j_r,k_r=0, 1, \ldots , n-1 \right).
\end{equation}
\end{theorem}
	
P\,r\,o\,o\,f.\, Necessity. Theorem~\ref{t*-1} implies that once  \eqref{e-8}
have place the relations  1) and 2) are satisfied
(in particular it follows from \eqref{e-*00}  that
 $\dfrac{\partial ^2 R\left(f,f' \right)}{\partial b_{n-1}^2}=
 a_0^{n-1}2\prod\limits_{i=3}^n g_i\ne 0$).
	
The mentioned conditions are also sufficient. Indeed.
Since  $R\left(f,f' \right)=0$ it follows that $f$ possesses a multiple root, say
 $z_1=z_2=w$. Herewith, taking into account  \eqref{e-*00}
 we conclude that inequality  $ 0\ne\dfrac{\partial ^2 R\left(f,f' \right)}{\partial b_{n-1}^2}=a_0^{n-1}2\prod\limits_{i=3}^n g_i$ implies
  $\prod\limits_{i=3}^n g_i\ne 0$, i.e. multiplicity of $w$
 is not greater than  2 and $f\left(z \right)$ does not possess any other multiple roots.

Relations  \eqref{e-*0*0} follow directly from equalities  \eqref{e-*00},
since in our case $s=2, k=1$. The proof is complete.

\medskip	

\emph{Remark}. \ 1)  Expressions  \eqref{e-*00} imply that once
a polynomial has the only one multiple root of multiplicity 2
 (all the other roots are simple) then all the second partial derivatives
 of the resultant  \eqref{e-3}
 are non-zero simultaneously.
 Thus one can verify the validity of relation
 \eqref{e-8} only for a one (any) of the second partial derivatives
 of the resultant  \eqref{e-*00}, for example, for~$\dfrac{\partial ^2 R}{\partial b_{n-1}^2}$.

2) 	 Once $R\left(f,f' \right)=0$  (i.e. there is at least one multiple root of $f$)
it follows that all the first partial derivatives are zero $\partial R\left(f,f' \right)/\partial b_j\ \left(j=0,1,\ldots , n-1 \right)$ (Theorem~\ref{t*-1}, equality \eqref{e-*0}).
Therefore verification of these expressions does not
enter necessary and sufficient conditions of satisfaction of relation~\eqref{e-8}.

\medskip

Theorem~\ref{t-6}	 is also a corollary of Lemma~\ref{l**-1}.
Namely the next refinement of Theorem~\ref{t-6} has place.

\begin{theorem}\label{t-6-n}
If a polynomial
 $f\left(z \right)=\sum\limits_{i=0}^n {a_i{z^{n-i}}}\  \left(a_0\neq 0 \right)$
 can be represented in the form
\begin{equation}\label{e-8-pp}
	f\left(z \right)=a_0{{\left(z-w \right)}^s}\left(z-{z_{s+1}} \right)\left(z-{z_{s+2}} \right)\ldots \left(z-{z_{n-1}} \right)\left(z-{z_n} \right),
\end{equation}
where  $2 \leq s< n;$ \   $w\neq z_i$ $($i.e.  $w$ is multiple root of multiplicity $s$$)$,

then the following relations are true$:$

$1)$ $R\left(f,f^{(k)} \right)=0, \ 1   \leq k < s;$

$2)$ for $R=R\left(f,f^{(s-1)} \right)$
the vector of derivatives  $\left[\frac{\partial R}{\partial {a_0}},\frac{\partial R}{\partial {a_1}},
\ldots , \frac{\partial R}{\partial {a_n}}\right]$ has the form
\begin{equation}\label{e-l**-1**p}
\left[\frac{\partial R}{\partial {a_0}},\frac{\partial R}{\partial {a_1}},
\ldots , \frac{\partial R}{\partial {a_n}}\right] =
\gamma (s,n) \prod\limits_{i=2}^{n-s+1} f_i\,\cdot \,[w^n,w^{n-1}, ... , w, 1].
\end{equation}
where  $\gamma (s,n)={{\left(-1 \right)}^{(n-s+1)n}}
\left(n(n-1)\cdot \ldots \cdot (n-s+2)a_0\right)^n$.

If   $ a_n\neq 0$  and all the у roots  $z_i, \ i= s+1, \ldots, n$
are of multiplicities less that  $s$, then all the derivatives
$\frac{\partial R}{\partial {a_j}}, \ j= 0,1,\ldots , n;$ are non-zero$;$
and therefore the relations \eqref{e-21} are true.
\end{theorem}

P\,r\,o\,o\,f.\, The argument in essence is the same as in the proof of
Proposition~\ref{p-3}.

A polynomial   $f$ possesses a multiple root of multiplicity $s$
iff  $R\left(f,f^{(k)} \right)=0, \ 1   \leq k < s$.  Let  $z_1=z_2 = \ldots =z_s=w$
be a multiple root of  $f$. This root  $w=y_1$ is a root of derivatives
 $f^{(k)}, \ 1   \leq k < s$.
In accordance with  \eqref{e-l**-1*} (for $k=s-1$) we have

   \begin{equation}\label{e-l**-1*a}
 \frac{\partial R}{\partial {a_j}} = {{\left(-1 \right)}^{(n-s+1)n}}b_0^n w^{n-j}  \prod\limits_{i=2}^{n-s+1} f_i, \ \ \ j=0, ... , n;
 \end{equation}
 where  $f_i$ are given by  \eqref{e-44} for  $k=s-1$.

Since  $b_0 =n(n-1)\cdot \ldots \cdot (n-s+2)a_0$ formula  \eqref{e-l**-1*a}
coincides with \eqref{e-l**-1**p}.

\smallskip

If
 $a_n\neq 0$ then all the roots of $f$ are non-zero  (in particular, $w\neq 0$).
 Herewith
  $b_0 =n(n-1)\cdot \ldots \cdot (n-s+2)a_0\neq 0$,
and  \eqref{e-l**-1*a} implies
\begin{equation*}
\frac{\partial R}{\partial {a_j}}  \neq 0 \ \ \Leftrightarrow \ \
\prod\limits_{i=2}^{n-s+1} f_i \neq 0.
\end{equation*}
Thus
\begin{equation}\label{e-t-a-0}
\frac{\partial R}{\partial {a_j}}  \neq 0 \ \ \Leftrightarrow \ \ f_i= f(y_i) \neq 0, \qquad i=2,\ldots , n-s+1.
\end{equation}
If  $z_i, \ i= s+1, \ldots , n$ are the roots of  $f$ of multiplicities  $ \leq s-1$
they are not the roots of the derivative  $f^{(s-1)}$.
Therefore, a common root of  $f$ and its derivative $f^{(s-1)}$
is only $y_1=w=z_1=z_2 = \ldots =z_s$. Hence
$f_i= f(y_i) \neq 0, \qquad i=2,\ldots , n-s+1$ and it follows that
 $\frac{\partial R}{\partial {a_j}}  \neq 0$ (in view of \eqref{e-t-a-0}).

Finally relations  \eqref{e-21} (taking into account \eqref{e-*})
follow directly from \eqref{e-l**-1**p}.
The proof is complete.

\bigskip

\emph{Remark}. \ 	Formula \eqref{e-l**-1**p} shows, in particular,
that a number of expressions of the form mentioned in  \eqref{e-21}
does not depend on $s$ (as in Theorem~\ref{t-6}) and is defined only by
degree of the initial polynomial.

\bigskip

The next generalization of Theorem~\ref{t-2} follows from Theorem~\ref{t*-1}
and is a certain alternative to Theorem~\ref{t-6-n} (Theorem~\ref{t-6}).


\begin{theorem}\label{t-8-n}
A polynomial
 $f\left(z \right)=\sum\limits_{i=0}^n{a_i {z^{n-i}}}\ \ \left(a_0\ne 0, \  a_n\ne 0\right)$
 can be represented in the form
\begin{equation}\label{e-56-n}
	f\left( z \right)={{a}_{0}}{{\left( z-w \right)}^{s}}\left( z-{{z}_{s+1}} \right)\left( z-{{z}_{s+2}} \right)\ldots \left( z-{{z}_{n-1}} \right)\left( z-{{z}_{n}} \right),
\end{equation}
where  $2\le s<n,\quad {{z}_{i}}\ne {{z}_{j}},\quad i\ne j$ and $w\ne {{z}_{i}};$
iff the following relations are satisfied$:$
	
	$1)$ $R\left( f,{{f}^{\left( s-1 \right)}} \right)=0$,
	
	\medskip
	$2)$ $R\left( f,{{f}^{\left( s \right)}} \right)\ne 0$,
	
	\medskip
	$3)$ $\dfrac{{{\partial }^{s}}R\left( f,{f}' \right)}{\partial b_{n-1}^{s}}\ne 0$,
where  ${{b}_{j}}\  \left( j=0, 1,\ldots , n-1 \right)$ are the coefficients of the
polynomial $g\left(z \right)=f'\left(z \right)$.

Herewith  $w$ satisfies the following relations$:$

\begin{equation}\label{e-60-n}
 	\frac{{{\partial }^{s}}R\left( f,{f}' \right)}{\partial {{b}_{{{j}_{s}}}}\ldots \partial {{b}_{{{j}_{1}}}}} :
\frac{{{\partial }^{s}}R\left( f,{f}' \right)}{\partial {{b}_{{{k}_{s}}}}\ldots \partial {{b}_{{{k}_{1}}}}} = {w}^{\left( {{k}_{s}}+\ldots +{{k}_{1}} \right)-\left( {{j}_{s}}+\ldots +{{j}_{1}} \right)}, \quad
  \left( {{j}_{r}, {k}_{r}}=0, 1, \ldots , n-1 \right).
\end{equation}
\end{theorem}
	
P\,r\,o\,o\,f.\, The necessity follows from Theorem~\ref{t*-1}.
Indeed, once the equality \eqref{e-56-n} takes place
relations  1) и 2) are true.
The equality  \eqref{e-*00} and Theorem~\ref{t*-1}  imply
	$$
	\frac{{{\partial }^{s}}R\left( f,{f}' \right)}{\partial b_{n-1}^{s}}=a_{0}^{n-1}s!\,\prod\limits_{i=s+1}^{n}{{{g}_{i}}}\ne 0.
	$$
	
Sufficiency. Since  $R\left( f,{{f}^{\left( s-1 \right)}} \right)=0$
it follows that  $f\left( z \right)$ possesses at least one multiple root of
multiplicity at least $s$, and inequality
$R\left( f,{{f}^{\left( s \right)}} \right)\ne 0$ implies that
$f\left( z \right)$  does not possess roots of multiplicity  $s+1$ or greater.
	
By  \eqref{e-*00} (Theorem~\ref{t*-1}) the partial derivatives have the form
\begin{equation}\label{e-57-n}
 	\frac{{{\partial }^{s}}R\left( f,{f}' \right)}{\partial {{b}_{{{j}_{s}}}}\ldots \partial {{b}_{{{j}_{1}}}}}=a_{0}^{n-1}s!\cdot {{w}^{s\left( n-1 \right)-\left( {{j}_{s}}+\ldots +{{j}_{1}} \right)}}\cdot \prod\limits_{i=s+1}^{n}{{{g}_{i}}}\quad \left( {{j}_{k}}=0, 1, \ldots , n-1 \right)
\end{equation}
(${{g}_{i}}\equiv g\left( {{z}_{i}} \right)$ are defined by  \eqref{e-4-}).
As $a_n\ne 0$ we have that $w\neq 0$, and therefore all the derivatives in
\eqref{e-57-n} are simultaneously  zero iff  $\prod\limits_{i=s+1}^{n}{{{g}_{i}}}=0$,
i.e. in the case when the original polynomial $f$ the original polynomial except
the root ${{z}_{1}}={{z}_{2}}=\ldots ={{z}_{s}}=w$ has another multiple root
of multiplicity at least~2. If $f\left( z \right)$ does not possess the second
 multiple root then all the partial derivatives
 $\dfrac{{{\partial }^{s}}R\left( f,{f}' \right)}{\partial {{b}_{{{j}_{s}}}}
 \ldots \partial {{b}_{{{j}_{1}}}}}$ are non-zero simultaneously.
 Therefore it suffice to verify any of them, for example,  $\dfrac{{{\partial }^{s}}R\left( f,{f}' \right)}{\partial b_{n-1}^{s}}$.

Relations  \eqref{e-60-n} follow directly from \eqref{e-57-n}.
The proof is complete.

\bigskip

\emph{Remark}. \ 	For the situation considered in the theorem relations \eqref{e-60-n}
are nothing else than relations \eqref{e-*00--}.

2) Once a polynomial $f$ possesses a root of multiplicity $s$ then for all
$1\le k<s$  the resultants  $R\left( f,{{f}^{\left( k \right)}} \right)=0$
automatically, therefore verification of these equalities does not enter
necessary and sufficient conditions of Theorem~\ref{t-8-n}.
	
	3) Condition 2) in Theorem~\ref{t-8-n} is redundant, since once $f\left( z \right)$
possesses a multiple root of multiplicity $s+1$ or greater  condition 3)
will not be satisfied (as ${{g}_{s+1}}=0$). Nevertheless often in practice
in the process of ``recognizing'' about the highest degree of multiplicity
this condition is  easier to verify than condition 3) (command for resultant calculation
does exist in mathematical pockets and is easily implemented in a  cycle of
iterating through the values of $R\left( f,{{f}^{\left( k \right)}} \right)$).
Once $s$ is determined it is advisable to check condition 3).

We will use this observation in examples presented in the next section.

\section{Examples}\label{s-3}

Here we give a number of  examples illustrating the obtaining of
concrete rational expressions for calculating multiple roots of
the greatest multiplicity based on the theorems from the previous section.

The formulas are rather cumbersome.
Therefore, we present their explicit form for polynomials of low degrees,
paying attention, first of all, to the qualitative difference
in the formulas derived on the basis of Theorems~\ref{t-5},~\ref{t-6} and
Theorems~\ref{t-2},~\ref{t-8-n}.

It is worth noting the high complexity of the manual calculation process, therefore, it is desirable to carry out all symbolic transformations in computer algebra systems.
We conducted them using the computer mathematics system \textit{Maple}.

\smallskip

Let us start from illustration of Proposition~\ref{p-3} (Theorem~\ref{t-5}).

Let a polynomial of the third degree
	$$
	f\left( z \right)=z^3+a_1z^2+a_2z+a_3\quad \left( a_0=1,\ a_3\ne 0 \right),
	$$
satisfy the conditions of Proposition~\ref{p-3}, i.e. can be represented in the form
$$
f\left(z \right)={{\left(z-w \right)}^2}\left(z-z_3\right).
$$
Then due to \eqref{e-l**-1**p-} the root $w$ satisfies the relations:

\begin{equation*}
\begin{split}
  w\quad &=
\quad\frac{\partial R}{\partial a_0}:\frac{\partial R}{\partial a_1}=-\frac{4a_1^3a_3-a_1^2a_2^2-9a_1a_2a_3+2a_2^3}{6a_1^2a_3-a_1a_2^2-9a_2a_3} \quad = \\
\quad &= \quad
\frac{\partial R}{\partial a_1}:\frac{\partial R}{\partial a_2}=-\frac{6a_1^2a_3-a_1a_2^2-9a_2a_3}{a_1^2a_2+9a_1a_3-6a_2^2}\quad = \\
\quad &= \quad
\frac{\partial R}{\partial a_2}:\frac{\partial R}{\partial a_3}=-\frac{a_1^2a_2+9a_1a_3-6a_2^2}{2a_1^3-9a_1a_2+27a_3},
\end{split}
\end{equation*}
where $R=R\left( f,f' \right)$.

\smallskip
	
Similarly for a polynomial of the fourth degree
$$
f\left( z \right)=z^4+a_1z^3+a_2z^2+a_3z+a_4\quad \left( a_0=1,\ a_4\ne 0 \right),
$$
possessing a root  $w$ of multiplicity 2 and other simple roots, i.e. having the form
$$
f\left(z \right)={{\left(z-w \right)}^2}\left(z-z_3 \right)\left(z-z_4 \right),
$$
the following relations are true
\begin{equation*}
\begin{split}
  w\quad &=
	\frac{\partial R}{\partial a_0}:\frac{\partial R}{\partial a_1}=-\left( 81a_1^4a_4^2-54a_1^3a_2a_3a_4+12a_1^3a_3^3+12a_1^2a_2^3a_4-3a_1^2a_2^2a_3^2-288a_1^2a_2a_4^2+ \right.\\
\quad & 	+12a_1^2a_3^2a_4+160a_1a_2^2a_3a_4-36a_1a_2a_3^3-32a_2^4a_4+8a_2^3a_3^2+192a_1a_3a_4^2+128a_2^2a_4^2-\\
\quad &
	-144a_2a_3^2a_4\left. +27a_3^4 \right)/\left( 2\left( 54a_1^3a_4^2- \right. \right.27a_1^2a_2a_3a_4+6a_1^2a_3^3+4a_1a_2^3a_4-a_1a_2^2a_3^2-\\
\quad &
	-144a_1a_2a_4^2+6a_1a_3^2a_4+40a_2^2a_3a_4-9a_2a_3^3\left. +\left. 96a_3a_4^2 \right) \right)\quad = \\
\quad &=
\frac{\partial R}{\partial a_1}:\frac{\partial R}{\partial a_2}=-\left( 54a_1^3a_4^2-27a_1^2a_2a_3a_4+6a_1^2a_3^3+4a_1a_2^3a_4-a_1a_2^2a_3^2-144a_1a_2a_4^2+ \right.\\
\quad &
	+\,6a_1a_3^2a_4+40a_2^2a_3a_4-9a_2a_3^3+\left. 96a_3a_4^2 \right)/\left( 9a_1^3a_3a_4-6a_1^2a_2^2a_4+a_1^2a_2a_3^2+72a_1^2a_4^2- \right.\\
\quad & -\,80a_1a_2a_3a_4+9a_1a_3^3+32a_2^3a_4-6a_2^2a_3^2-128{a_2}a_4^2\left. +72a_3^2a_4 \right)
\quad = \\
\quad &=
\frac{\partial R}{\partial a_2}:\frac{\partial R}{\partial a_3}=\left( 9a_1^3a_3a_4-6a_1^2a_2^2a_4+a_1^2a_2a_3^2+72a_1^2a_4^2-80a_1a_2a_3a_4 \right.+9a_1a_3^3+32a_2^3a_4-\\
\quad &
	-\,6a_2^2a_3^2-128a_2a_4^2+\left. 72a_3^2a_4 \right)/\left( 9a_1^3a_2a_4-6a_1^3a_3^2+a_1^2a_2^2a_3-6a_1^2a_3a_4-40a_1a_2^2a_4+ \right.\\
\quad &
	+27a_1a_2a_3^2-4a_2^3{a_3}-96{a_1}a_4^2+144a_2a_3a_4-\left. 54a_3^3 \right)\quad = \\
\quad &= \frac{\partial R}{\partial a_3}:\frac{\partial R}{\partial a_4}=-\left( 9a_1^3a_2a_4-6a_1^3a_3^2+a_1^2a_2^2a_3-6a_1^2a_3a_4-40a_1a_2^2a_4+27a_1a_2a_3^2 \right.-\,4a_2^3a_3-\\
\quad &
	-\,96a_1a_4^2+144a_2a_3a_4-\left. 54a_3^3 \right)/\left( +27a_1^4a_4-9a_1^3a_2a_3+2a_1^2a_2^3-144a_1^2a_2a_4+ \right.\\
\quad &
	+3a_1^2a_3^2+40a_1a_2^2a_3-8a_2^4+192a_1a_3a_4+128a_2^2a_4-72a_2a_3^2-\left. 384a_4^2 \right).
\end{split}
\end{equation*}

	\bigskip

Let us now illustrate Theorem~\ref{t-2} using the same examples.

Let a polynomial of the third degree
	$$
	f\left( z \right)=z^3+a_1z^2+a_2z+a_3\quad \left( a_0=1,\ a_3\ne 0 \right),
	$$
satisfy the conditions of Theorem~\ref{t-2}, i.e. has the form
$$
f\left(z \right)={{\left(z-w \right)}^2}\left(z-z_3\right).
$$
Since (see \eqref{e-*0*0})
\begin{equation}\label{e-*0*0-}
 \frac{\partial ^2 R}{\partial b_{j_2}\partial b_{j_1}} :
 \frac{\partial ^2 R}{\partial b_{k_2}\partial b_{k_1}} =
 w^{(k_2+k_1) -(j_2 +j_1)}, \qquad
  \left(j_r,k_r=0, 1, \ldots , n-1 \right),
\end{equation}
for $(k_2+k_1) -(j_2 +j_1) = 1$ we obtain the explicit formulas for  $w$.

Below are all the formulas for the root $w$
that follow from here.

\begin{equation*}
\begin{split}
  w\quad &= \quad	\frac{\partial ^2 R}{\partial b_0^2}:\frac{\partial ^2 R}{\partial b_0\partial b_1}=\frac{2a_2^3-8a_1 a_2 a_3 +18a_3^2}{4a_1^2 a_3- a_1 a_2^2-3a_2 a_3}\quad = \\
& =\quad 	\frac{\partial ^2 R}{\partial b_0\partial b_1}:\frac{\partial ^2 R}{\partial b_0\partial b_2}=\frac{4a_1^2 a_3- a_1 a_2^2-3a_2 a_3}{2a_2^2-6a_1 a_3}\quad = \\
&	= \quad \frac{\partial ^2 R}{\partial b_0\partial b_1}:\frac{\partial ^2 R}{\partial b_1^2}=\frac{4a_1^2 a_3- a_1 a_2^2-3a_2 a_3}{2a_2^2-6a_1 a_3}\quad = \\
&= \quad 	\frac{\partial ^2 R}{\partial b_1^2}:\frac{\partial ^2 R}{\partial 	{b_1}\partial b_2}=\frac{2a_2^2-6a_1 a_3}{9{a_3}-a_1 a_2}\quad = \\
&	= \quad \frac{\partial ^2 R}{\partial b_0\partial b_2}:\frac{\partial ^2 R}{\partial b_1\partial b_2}=\frac{2a_2^2-6a_1 a_3}{9 a_3 -a_1 a_2}\quad = \\
& =\quad 	\frac{\partial ^2 R}{\partial b_1\partial b_2}:\frac{\partial ^2 R}{\partial b_2^2}=\frac{9 a_3-a_1 a_2}{2a_1^2-6 a_2}.
\end{split}
\end{equation*}

One observes  that some of the expressions obtained coincide.
To make sure that they are valid in general, it is enough to substitute in them the Viete's
formulas  written taking into account the multiple root, that is:
$a_1=-2z_1-z_3,\ a_2=z_1^2+2 z_1 z_3,\ a_3=-z_1^2z_3$.
 Then, after simplifications, each fraction will turn into $z_1$.

Note also that the formulas presented above (as well as those obtained further
 in  the article) are significantly simplified if $a_1=0$.
  We can always reduce the consideration to this situation.
  Indeed, using the well-known substitution $z=y-a_1/n$ the algebraic equation
 ${z^n}+{a_1}{z^{n-1}}+a_2 {z^{n-2}}+\ldots +a_n=0$  reduces to the equation
 ${y^n}+c_2{z^{n-2}}+c_3{z^{n-3}}+\ldots +{c_{n-1}}z+{c_n}=0$.

Now let us illustrate Theorem~\ref{t-2} by example of the fourth degree polynomial
$$
f\left( z \right)=z^4+a_1z^3+a_2z^2+a_3z+a_4\quad \left( a_0=1,\ a_4\ne 0 \right),
$$
possessing a root  $w$ of multiplicity 2 and simple other roots, i.e. having the form
$$
f\left(z \right)={{\left(z-w \right)}^2}\left(z-z_3 \right)\left(z-z_4 \right).
$$
Below we present some  formulas for calculating the root $w$
based on formulas \eqref{e-*0*0-}
 for  $(k_2+k_1) -(j_2 +j_1) = 1$.

\begin{equation*}
\begin{split}
  w\quad &= \quad
		\frac{\partial ^2 R}{\partial b_0^2}:\frac{\partial ^2 R}{\partial b_0\partial b_1}=\left(54a_1^2 a_2 a_4^2 \right.-36 a_1 a_2^2a_3 a_4+8a_1 a_2a_3^3+8a_2^4 a_4 -2a_2^3a_3^2-120a_1 a_3a_4^2- \\
	& \quad
	-80\,a_2^2a_4^2+94 a_2 a_3^2 a_4 -18\,a_3^4+\left. 192a_4^3 \right)/\left(-27a_1^3a_4^2+18a_1^2a_2 a_3 a_4-4a_1^2a_3^3-4 a_1 a_2^3 a_4 + \right. \\
	& \quad + a_1 a_2^2a_3^2+48\,a_1 a_2a_4^2-7 a_1 a_3^2 a_4 -12a_2^2a_3 a_4+3 a_2 a_3^3\left. -16 a_3 a_4^2 \right) \quad = \\
	 &= \quad
	\frac{\partial ^2 R}{\partial b_1^2}:\frac{\partial ^2 R}{\partial b_1\partial b_2}=\frac{36a_1^2a_4^2-28a_1 a_2a_3 a_4+6 a_1 a_3^3+8a_2^3 a_4 -2a_2^2a_3^2-32 a_2 a_4^2+12a_3^2a_4 }{3a_1^2a_3 a_4-4 a_1 a_2^2 a_4 +a_1 a_2a_3^2-48 a_1 a_4^2+32a_2 a_3 a_4-9a_3^3}\quad = \\
 &= \quad
	\frac{\partial ^2 R}{\partial b_1\partial b_2}:\frac{\partial ^2 R}{\partial b_2^2}=\frac{3a_1^2a_3 a_4-4 a_1 a_2^2 a_4 +a_1 a_2a_3^2-48 a_1 a_4^2+32a_2 a_3 a_4 -9a_3^3}{6a_1^2a_2 a_4-2a_1^2a_3^2-8a_1 a_3 a_4 -16a_2^2 a_4 +6 a_2 a_3^2+64a_4^2}\quad = \\
&= \quad 	
	\frac{\partial ^2 R}{\partial b_1\partial b_3}:\frac{\partial ^2 R}{\partial b_2\partial b_3}=\frac{6a_1^2a_2 a_4-2a_1^2a_3^2-8\,a_1 a_3 a_4 -16a_2^2 a_4 +6 a_2 a_3^2+64a_4^2}{a_1^2 a_2 a_3-9a_1^3 a_4 +32a_1 a_2 a_4 +3 a_1 a_3^2-4a_2^2 a_3 -48a_3 a_4} \quad = \\
 &= \quad
	\frac{\partial ^2 R}{\partial b_2\partial b_3}:\frac{\partial ^2 R}{\partial b_3^2}=\frac{a_1^2a_2 a_3-9a_1^3 a_4+32a_1 a_2 a_4 +3 a_1 a_3^2-4a_2^2 a_3 -48a_3 a_4}{6a_1^3 a_3 -2a_1^2a_2^2+12a_1^2 a_4-28a_1 a_2 a_3 +8a_2^3-32a_2 a_4+36a_3^2} \ .
\end{split}
\end{equation*}

Now let us illustrate Theorem~\ref{t-8-n} by example of the fourth degree polynomial
$$
f\left( z \right)=z^4+a_1z^3+a_2z^2+a_3z+a_4\quad \left( a_0=1,\ a_4\ne 0 \right),
$$
possessing a root  $w$ of multiplicity 3 and a simple root, i.e.
having the form
$$
f\left(z \right)={{\left(z-w \right)}^3}\left(z-z_4 \right).
$$
In this situation formula  \eqref{e-60-n} reduces to
\begin{equation}\label{e-60-n-}
 	\frac{{{\partial }^{3}}R\left( f,{f}' \right)}{\partial {{b}_{{{j}_{3}}}}\partial {{b}_{{{j}_{2}}}} \partial {{b}_{{{j}_{1}}}}} :
\frac{{{\partial }^{s}}R\left( f,{f}' \right)}{\partial {{b}_{{{k}_{3}}}}\partial {{b}_{{{k}_{2}}}} \partial {{b}_{{{k}_{1}}}}} = {w}^{\left( {{k}_{3}}+{{k}_{2}} +{{k}_{1}} \right)-\left( {{j}_{3}}+{{j}_{2}} +{{j}_{1}} \right)}, \quad
  \left( {{j}_{r}, {k}_{r}}=\, 0, 1, 2, 3 \right).
\end{equation}
And for $\left( {{k}_{3}}+{{k}_{2}} +{{k}_{1}} \right)-\left( {{j}_{3}}+{{j}_{2}} +{{j}_{1}} \right) = 1$ one obtains explicit formulas for
 $w$.

Thus, even on a simple example of a polynomial of the fourth degree
one can see a variety of explicit rational formulas for a multiple root (more than twenty different symbolic expressions). Here are some of them.

\begin{equation*}
\begin{split}
  w\quad &= \quad
		\frac{\partial ^3 R}{\partial b_{0}^3}:\frac{\partial ^3 R}{\partial b_0^2\partial b_1}=\frac{3\left(5 a_2 a_3^2 a_4-6a_1 a_3a_4^2-4a_2^2a_4^2-a_3^4+16a_4^3 \right)}{12a_1 a_2a_4^2-a_1 a_3^2 a_4 -4a_2^2a_3 a_4+ a_2 a_3^3-8 a_3 a_4^2} \quad = \\	
&= \quad
	\frac{\partial ^3 R}{\partial b_0^2\partial b_1}:\frac{\partial ^3 R}{\partial b_0^2\partial b_2}=\frac{12a_1 a_2a_4^2-a_1 a_3^2 a_4 -4 a_2^2a_3 a_4+a_2a_3^3-8 a_3 a_4^2}{9a_1^2a_4^2-10a_1 a_2a_3 a_4+2 a_1 a_3^3+4a_2^3 a_4-a_2^2a_3^2-16 a_2 a_4^2+7a_3^2 a_4} \quad = \\	
&= \quad
	\frac{\partial ^3 R}{\partial b_0^2\partial b_2}:\frac{\partial ^3 R}{\partial b_0^2\partial b_3}=\frac{9a_1^2a_4^2-10a_1 a_2a_3 a_4+2 a_1 a_3^3+4a_2^3{a_4}-a_2^2a_3^2-16 a_2 a_4^2+7a_3^2 a_4}{3a_1^2a_3 a_4-4 a_1 a_2^2 a_4 +a_1 a_2a_3^2-24 a_1 a_4^2+20a_2 a_3 a_4-6a_3^3} \quad = \\
&= \quad
	\frac{\partial ^3 R}{\partial b_0^2\partial b_1}:\frac{\partial ^3 R}{\partial b_0^2\partial b_3}=\frac{12a_1 a_2a_4^2- a_1 a_3^2 a_4 -4a_2^2a_3 a_4+ a_2 a_3^3-8 a_3 a_4^2}{4a_1 a_2a_3 a_4-9a_1^2a_4^2- a_1 a_3^3+a_3^2a_4} \quad = \\
&= \quad
	\frac{\partial ^3 R}{\partial b_0\partial b_1\partial b_2}:\frac{\partial ^3 R}{\partial b_0\partial b_1\partial b_3}=\frac{\partial ^3 R}{\partial b_0\partial b_1\partial b_2}:\frac{\partial ^3 R}{\partial b_0\partial b_2^2}\quad = \\
&= \quad \frac{3a_1^2a_3 a_4-4 a_1 a_2^2 a_4 +a_1 a_2a_3^2+8a_2 a_3 a_4-3a_3^3}{2\left(3a_1^2a_2 a_4-a_1^2a_3^2-6a_1 a_3 a_4-4a_2^2 a_4 +2 a_2 a_3^2+16a_4^2 \right)}\quad = \\
&= \quad
	\frac{\partial ^3 R}{\partial b_2\partial b_3^2}:\frac{\partial ^3 R}{\partial b_3^3}=\frac{a_1^2 a_2 +2a_1 a_3-4\,a_2^2+16 a_4 }{3\left(4a_1 a_2-a_1^3-8 a_3 \right)}\quad = \\
&= \quad
	\frac{\partial ^3 R}{\partial b_1^2\partial b_2}:\frac{\partial ^3 R}{\partial b_1^2\partial b_3}=\frac{\partial ^3 R}{\partial b_1^2\partial b_2}:\frac{\partial ^3 R}{\partial b_1\partial b_2^2}=\frac{4a_2^2 a_4 -2a_1 a_3 a_4-a_2 a_3^2-16a_4^2}{ a_1 a_3^2-4a_1 a_2 a_4 +8a_3 a_4}
\quad = \\
&= \quad
	\frac{\partial ^3 R}{\partial b_1\partial b_2^2}:\frac{\partial ^3 R}{\partial b_1\partial b_2\partial b_3}=\frac{\partial ^3 R}{\partial b_1^2\partial b_3}:\frac{\partial ^3 R}{\partial b_1\partial b_2\partial b_3}=\frac{\partial ^3 R}{\partial b_1\partial b_2^2}:\frac{\partial ^3 R}{\partial b_2^3}=\frac{a_1 a_3^2-4a_1 a_2 a_4 +8a_3 a_4}{3\left(a_1^2 a_4-a_3^2 \right)}
\quad = \\
&= \quad
	\frac{\partial ^3 R}{\partial b_1\partial b_2\partial b_3}:\frac{\partial ^3 R}{\partial b_1\partial b_3^2}=\frac{\partial ^3 R}{\partial b_2^3}:\frac{\partial ^3 R}{\partial b_2^2\partial b_3}=\frac{\partial ^3 R}{\partial b_1\partial b_2\partial b_3}:\frac{\partial ^3 R}{\partial b_2^2\partial b_3}=\frac{3\left(a_1^2 a_4-a_3^2 \right)}{4a_2 a_3-a_1^2 a_3-8a_1 a_4}
\quad = \\
&= \quad
	\frac{\partial ^3 R}{\partial b_0\partial b_1^2}:\frac{\partial ^3 R}{\partial b_1^3}=\frac{4a_1 a_2a_3 a_4-9a_1^2a_4^2-a_1a_3^3+a_3^2a_4}{3\left(8a_1a_4^2-4a_2 a_3 a_4+a_3^3 \right)}.
\end{split}
\end{equation*}
	
Under substitution of Viete's relations
$$
a_1=-3z_1-z_4, \ a_2=3z_1^2+3z_1z_4, a_3=-z_1^3-3 z_1^2 z_4,\ a_4=z_1^3 z_4
$$
into all the listed fractions	by means of math package
we make sure that after the reduction each fraction will turn into $z_1$.

\end{document}